\documentclass[a4paper,12pt,reqno]{amsart}
\usepackage{amsmath}
\usepackage{amssymb}
\usepackage{amsthm}

\usepackage{amscd}
\usepackage{amsfonts}
\usepackage{setspace}
\usepackage{version}

\usepackage{hyperref}

\title{Better than squareroot cancellation in number theory}
\author{Adam J Harper}
\address{Mathematics Institute, Zeeman Building, University of Warwick, Coventry CV4 7AL, England}
\email{A.Harper@warwick.ac.uk}
\date{29th December 2025}
\thanks{This research was funded in part by the Engineering and Physical Sciences Research Council of the United Kingdom [grant EP/V055755/1]. For the purpose of open access, the author has applied a Creative Commons Attribution (CC-BY) licence to any Author Accepted Manuscript version arising from this submission.}

\numberwithin{equation}{section}
\theoremstyle{plain}

\onehalfspace
\newcommand{\N}{\mathbb{N}}
\newcommand{\R}{\mathbb{R}}
\newcommand{\E}{\mathbb{E}}
\newcommand{\p}{\mathbb{P}}

\newcommand{\C}{\mathbb{C}}

\newcommand{\Echar}{\mathbb{E}^{\text{char}}}

\newtheorem{theorem}{Theorem}
\numberwithin{theorem}{section}

\numberwithin{proposition}{section}

\numberwithin{corollary}{section}

\newtheorem{conjecture}{Conjecture}
\numberwithin{conjecture}{section}

\begin{document}

\maketitle

\begin{abstract}
We give a short survey of the phenomenon of better than squareroot cancellation, specifically as it applies to averages of multiplicative character sums (such as $\frac{1}{r-1} \sum_{\chi \; \text{mod} \; r} |\sum_{n \leq x} \chi(n)|^{2q}$) thanks to their connection with so-called multiplicative chaos.

We focus on the number theoretic aspects of the arguments, and also touch on some possible applications.
\end{abstract}

\section{Introduction}

\subsection{Squareroot cancellation}
It is a familiar heuristic across much of mathematics that if one sums a sequence of oscillating terms, then ``typically'' the sum should have size ``around squareroot'' of the ``trivial bound''. This phenomenon often goes under the name of {\em squareroot cancellation}. Each of the highlighted parts of the sentence deserves some attention, but lots of cases should spring to mind that certainly conform to it.

The most fundamental example perhaps arises from taking a sequence $(\epsilon_n)_{n=1}^{\infty}$ of independent, fair $\pm 1$ coin tosses (so $\p(\epsilon_n = 1) = \p(\epsilon_n = -1) = 1/2$, independently for all $n$). Then a trivial bound for $\sum_{n \leq x} \epsilon_n$ is $\sum_{n \leq x} |\epsilon_n| = \lfloor x \rfloor$. But this is a sum of independent, mean zero, variance one random variables, and the most classical version of the Central Limit Theorem asserts that $\frac{1}{\sqrt{x}} \sum_{n \leq x} \epsilon_n$ converges in distribution to a standard normal random variable as $x \rightarrow \infty$. In particular, this means that for any given large $x$, we will have $c \leq \frac{1}{\sqrt{x}} |\sum_{n \leq x} \epsilon_n| \leq C$ with probability at least 0.99 (say), for appropriate constants $0 < c < C$. In this example, ``typically'' meant ``for most realisations'' of the random sequence $\epsilon_n$. Let us also note carefully that we have squareroot cancellation, but {\em not more than that}. The sum is typically no more than a large constant times $\sqrt{x}$ (the squareroot of the trivial bound of $\lfloor x \rfloor \approx x$), but it is also at least a small constant times $\sqrt{x}$. Another (perhaps better) way of thinking about this is that the variance of the random sum, namely $\E |\sum_{n \leq x} \epsilon_n|^2 = \sum_{n \leq x} \E \epsilon_n^2 = \sum_{n \leq x} 1 = \lfloor x \rfloor$, turns out to accurately reflect the typical behaviour of the sum: the variance suggests that we might usually have $|\sum_{n \leq x} \epsilon_n| \approx \sqrt{x}$, and the Central Limit Theorem confirms this is indeed the case.

\vspace{12pt}
In many settings, including those we shall discuss in more detail in this survey, the ``typical'' size of a sum has a similar meaning as in the above example. That is, we actually have a family of many different but related sums, and we want to obtain size bounds that hold for most sums in the family. We emphasise, however, that the kinds of families we shall focus on will be deterministic: they will come from number theory and involve no ``genuine'' randomness, although it will turn out that their behaviour is closely tied up with probabilistic ideas.

In other cases, one has a single specific sum, and expects it by itself to exhibit cancellation. Or one might have a family of sums, but expect all of them to exhibit cancellation unless there is some ``obvious'', checkable reason why that cannot hold. For example, we can consider the sum of the M\"obius function $\mu(n)$, an object that we shall discuss extensively. Recall that $\mu : \N \rightarrow \R$ is defined by $\mu(n) = 0$ if $n$ has any repeated prime factors, and otherwise (i.e. if $n$ is squarefree) $\mu(n)$ is $\pm 1$ according to the parity of the number of prime factors of $n$. A trivial bound for $\sum_{n \leq x} \mu(n)$ is again $\sum_{n \leq x} |\mu(n)| \leq \sum_{n \leq x} 1 = \lfloor x \rfloor$ (and a slightly less trivial bound, taking account of the density of the squarefree numbers, would be $\approx (6/\pi^2)x$). It is known that $\sum_{n \leq x} \mu(n) = o(x)$ as $x \rightarrow \infty$, a deep result equivalent to the Prime Number Theorem. And it is conjectured, and equivalent to the Riemann Hypothesis, that one should have $|\sum_{n \leq x} \mu(n)| \leq x^{1/2 + o(1)}$ as $x \rightarrow \infty$. Since the M\"obius function is a specific deterministic function, we must be careful when trying even to guess its exact behaviour (a point we shall return to later). Using its connection with the Riemann zeta zeros and random matrix theory, Gonek has conjectured (see Ng's paper~\cite{ng}) that the largest fluctuations of $\sum_{n \leq x} \mu(n)$, as $x$ varies, should have order $\sqrt{x} (\log\log\log x)^{5/4}$.

\vspace{12pt}
Before moving on to our main objects of study, let us comment a bit more on the general meaning of ``around squareroot of the trivial bound''. As in the above discussion, it is clearly a subjective matter to identify ``the'' trivial bound in any given problem. But motivated by the probabilistic situation, if $(a_n)_{n=1}^{\infty}$ is a sequence of summands then (by the Cauchy--Schwarz inequality) we certainly always have $|\sum_{n \leq x} a_n| \leq \sqrt{\sum_{n \leq x} 1} \sqrt{\sum_{n \leq x} |a_n|^2} \leq \sqrt{x} \sqrt{\sum_{n \leq x} |a_n|^2}$, and if the size of the sum is actually $\approx \sqrt{\sum_{n \leq x} |a_n|^2}$ (i.e. saving a factor around $\sqrt{x}$) then we might consider ourselves to have squareroot cancellation. Notice that $\sum_{n \leq x} a_n \epsilon_n$ would again be a sum of independent random variables $a_n \epsilon_n$, with mean zero but now each having variance $|a_n|^2$. Under mild conditions (that no small collection of summands is too highly weighted relative to the others, see e.g. Chapter 7.2 of Gut~\cite{gut}), one again has a Central Limit Theorem {\em after renormalising by the standard deviation} $\sqrt{\sum_{n \leq x} |a_n|^2}$, and so $|\sum_{n \leq x} a_n \epsilon_n |$ will lie between $c\sqrt{\sum_{n \leq x} |a_n|^2}$ and $C\sqrt{\sum_{n \leq x} |a_n|^2}$ with probability 0.99. 

Let us finally observe that there is much overkill in invoking the Central Limit Theorem in this discussion. For upper bounds, Chebychev's inequality alone implies that
$$ \p\Biggl(|\sum_{n \leq x} a_n \epsilon_n | \geq \lambda \sqrt{\sum_{n \leq x} |a_n|^2} \Biggr) \leq \frac{\E|\sum_{n \leq x} a_n \epsilon_n |^2}{\lambda^2 \sum_{n \leq x} |a_n|^2} = \frac{1}{\lambda^2} \;\;\;\;\; \forall \lambda > 0 . $$
Thus $|\sum_{n \leq x} a_n \epsilon_n |$ must typically be at most $C\sqrt{\sum_{n \leq x} |a_n|^2}$ (but the Central Limit Theorem, or other large deviation results, would give stronger upper bounds for the probability). {\em Lower bounds require more than just the second moment}, but (importantly using the independence of the $\epsilon_n$) one can compute e.g. that
$$ \E|\sum_{n \leq x} a_n \epsilon_n |^4 = \sum_{\substack{n_1, n_2, n_3, n_4 \leq x, \\ n_1 = n_2 \; \text{and} \; n_3 = n_4, \\ \text{or} \; n_1 = n_3 \; \text{and} \; n_2 = n_4, \\ \text{or} \; n_1 = n_4 \; \text{and} \; n_2 = n_3}} a_{n_1} a_{n_2} \overline{a_{n_3} a_{n_4}} \leq 3 \left( \sum_{n \leq x} |a_n|^2 \right)^2 . $$
Combining this fourth moment estimate with the fact that $\E|\sum_{n \leq x} a_n \epsilon_n |^2 = \sum_{n \leq x} |a_n|^2$, and the Cauchy--Schwarz inequality, we get e.g. that
\begin{eqnarray}
\p\Biggl(|\sum_{n \leq x} a_n \epsilon_n | \geq \frac{1}{\sqrt{2}} \sqrt{\sum_{n \leq x} |a_n|^2} \Biggr) & \geq & \frac{(\E \textbf{1}_{|\sum_{n \leq x} a_n \epsilon_n | \geq \frac{1}{\sqrt{2}} \sqrt{\sum_{n \leq x} |a_n|^2}} |\sum_{n \leq x} a_n \epsilon_n |^2 )^2}{\E|\sum_{n \leq x} a_n \epsilon_n |^4} \nonumber \\
& \geq & \frac{(\frac{1}{2} \sum_{n \leq x} |a_n|^2 )^2}{3 (\sum_{n \leq x} |a_n|^2)^2} \geq \frac{1}{12} . \nonumber
\end{eqnarray}

\subsection{Better than squareroot cancellation}
In many number theoretic problems, proving anything close to squareroot cancellation (especially for individual sums) seems tremendously difficult, and even far weaker bounds constitute huge breakthroughs. For example, the best known unconditional estimate for $\sum_{n \leq x} \mu(n)$, equivalent to the very deep Vinogradov--Korobov zero-free region for the Riemann zeta function, is of the shape $x^{1-o(1)}$ rather than the conjectured $x^{1/2 + o(1)}$. In other problems, one can prove squareroot cancellation and knows this to be unimprovable. Thus the celebrated Weil bound~\cite{weil} for Kloosterman sums (and more general types of character sum) implies that $|\sum_{x=1}^{p-1} \exp(2\pi i \frac{ax + b\overline{x}}{p})| \leq 2\sqrt{p}$, if $p$ is a prime not dividing $ab$, and $\overline{x}$ denotes the inverse of $x$ modulo $p$. And it is known e.g. that $|\sum_{x=1}^{p-1} \exp(2\pi i \frac{x + b\overline{x}}{p})| \gg \sqrt{p}$ for a proportion 0.99 of values $1 \leq b \leq p-1$ (where $p$ is large), and in fact one understands the distribution of $\frac{1}{2\sqrt{p}} |\sum_{x=1}^{p-1} \exp(2\pi i \frac{x + b\overline{x}}{p})|$ as $b$ varies (thanks to the work of Katz~\cite{katz} on the vertical Sato--Tate law).

{\em In this survey we shall discuss a natural class of problems where it has been discovered, in recent years, that one actually obtains better than (i.e. more than) squareroot cancellation.} The extra saving one obtains is quite small--- perhaps indicating the subtlety of the underlying behaviour--- but has potentially very appealing consequences. Our goal here is to explain this, as well as some ideas from the proof of better than squareroot cancellation, and possible further developments in this programme of work.

\vspace{12pt}
We consider two families of sums: the sums $\sum_{n \leq x} a_n \chi(n)$, which one wants to understand for most Dirichlet characters $\chi$ modulo $r$ (a large prime, say); and the sums $\sum_{n \leq x} a_n n^{-it}$, which one wants to understand for most values of $0 \leq t \leq T$ (or possibly, for technical convenience, with some smooth weighting applied to the $t$ values).

Dirichlet characters $\chi(n)$, and the ``continuous characters'' $n \mapsto n^{-it}$, share the crucial property of {\em (total) multiplicativity}:
$$ \chi(nm) = \chi(n)\chi(m) \;\;\; \forall \; n,m , \;\;\;\;\; (nm)^{-it} = n^{-it} m^{-it} \;\;\; \forall \; n,m . $$
This property is very natural in the context of multiplicative number theory problems, such as the distribution of primes or smooth numbers, and many classical number theoretic functions are multiplicative. It turns out that multiplicativity can be a key driver of the better than squareroot cancellation phenomenon, so we will be interested in coefficients $a_n$ that are themselves multiplicative (or somewhat close to being so) as a function of $n$. With truly arbitrary coefficients (e.g. independent random coefficients) the phenomenon would disappear, as in our initial example of $\sum_{n \leq x} \epsilon_n$, which typically has squareroot cancellation and not more than that.

Dirichlet characters and continuous characters also both enjoy {\em orthogonality} properties (meaning, in particular, that their values are typically oscillatory). For Dirichlet characters modulo a prime $r$ the property is very clean: for any $n,m \in \N$ that are not divisible by $r$, we have
$$ \frac{1}{r-1} \sum_{\chi \; \text{mod} \; r} \chi(n) \overline{\chi(m)} = \textbf{1}_{n \equiv m \; \text{mod} \; r} , $$
where $\textbf{1}$ denotes the indicator function, and the sum is over the $r-1$ different characters. In particular, if we constrain $n,m$ to satisfy $1 \leq n,m < r$, {\em but crucially only with such a constraint}, then in fact
\begin{equation}\label{eqncharorth}
\frac{1}{r-1} \sum_{\chi \; \text{mod} \; r} \chi(n) \overline{\chi(m)} = \textbf{1}_{n = m} .
\end{equation}
When averaging over $0 \leq t \leq T$, we get approximate orthogonality:
$$ \frac{1}{T} \int_{0}^{T} n^{-it} m^{it} dt = \frac{1}{T} \int_{0}^{T} e^{-it\log(n/m)} dt = \textbf{1}_{n=m} + O\left(\frac{\textbf{1}_{n \neq m}}{T|\log(n/m)|} \right) . $$
Notice that if $1 \leq m < n$, say, then Taylor expansion of the logarithm implies the ``big Oh'' term can be bounded by $O(\frac{n}{T})$, so it will be small provided $m,n$ are rather smaller than $T$. With a smoothly weighted integral over $t$, one could improve this further.

A third important property to note is some form of {\em periodicity}, and resulting {\em symmetry} in certain sums involving $\chi(n)$ or $n^{-it}$. Again, this is cleanest in the Dirichlet character case. By definition, one has $\chi(n+r) = \chi(n)$ for all $n \in \N$ and all characters $\chi$ modulo $r$. Thus $\chi$ is sensitive to additive structure modulo $r$, or (the reverse viewpoint) $\chi$ is really a character on the finite multiplicative group of invertible residue classes mod $r$, with its definition extended to be a function on $\N$. All of this means that the additive Fourier coefficients of $\chi$ mod $r$ have very nice structure, which manifests itself in properties like the functional equations of the corresponding Dirichlet $L$-functions $L(s,\chi)$. For our purposes, a key consequence is the so-called ``Fourier flip'' (see e.g. section 10 of Granville and Soundararajan~\cite{gransoundlcs}): roughly speaking, and in a suitably averaged sense over $\chi$ and/or $x$, one has
$$ \frac{1}{\sqrt{x}} |\sum_{n \leq x} \chi(n)| \approx \sqrt{\frac{x}{r}} |\sum_{n \leq r/x} \chi(n)| , \;\;\;\;\; x \leq r . $$
So when the weights $a_n$ are identically 1, say, a long character sum of length $x \approx r$ is essentially the same, up to scaling, as a very short character sum of length $r/x \approx 1$. For continuous characters $n \mapsto n^{-it}$ there is no such perfect periodicity, but one obtains a similar ``flip'' corresponding to the functional equation of the Riemann zeta function: for large $t$, and in a suitably averaged sense, we get
$$ \frac{1}{\sqrt{x}} |\sum_{n \leq x} n^{-it}| \approx \sqrt{\frac{x}{|t|}} |\sum_{n \leq |t|/x} n^{-it}| , \;\;\;\;\; x \leq |t| . $$

\vspace{12pt}
With the above preamble, we can state some concrete results. Because of orthogonality, it is easy to see that for any $x < r$ we have $\frac{1}{r-1} \sum_{\chi \; \text{mod} \; r} |\sum_{n \leq x} \chi(n)|^2 = \lfloor x \rfloor$, as in the random case with $\epsilon_n$. With a little work, approximate orthogonality (or in fact the Montgomery--Vaughan mean value theorem, see e.g. Chapter 7 of Montgomery~\cite{mont}) similarly implies that
$$ \frac{1}{T} \int_{0}^{T} |\sum_{n \leq x} n^{-it}|^2 dt = \lfloor x \rfloor + O(x^{2}/T) . $$
The squareroot cancellation heuristic would suggest that for most $\chi$ and $t$ we should have $|\sum_{n \leq x} \chi(n)|, |\sum_{n \leq x} n^{-it}| \asymp \sqrt{x}$, for $x < r$ and $x < T$ respectively. But in fact, we get:
\begin{theorem}[Harper~\cite{harpertypicalchar}, 2023]\label{thmplainlow}
In the above setting with $r$ a large prime, uniformly for any $1 \leq x \leq r$ and any $0 \leq q \leq 1$ we have
$$ \frac{1}{r-1} \sum_{\chi \; \text{mod} \; r} |\sum_{n \leq x} \chi(n)|^{2q} \ll \Biggl(\frac{x}{1 + (1-q)\sqrt{\log\log(10L)}} \Biggr)^q , $$
where $L = L_r := \min\{x,r/x\}$.

Uniformly for any $1 \leq x \leq T$ (with $T \in \R$ large) and any $0 \leq q \leq 1$, we have
$$ \frac{1}{T} \int_{0}^{T} |\sum_{n \leq x} n^{it}|^{2q} dt \ll \Biggl(\frac{x}{1 + (1-q)\sqrt{\log\log(10L_T)}} \Biggr)^q  , $$
where $L_T := \min\{x,T/x\}$.
\end{theorem}

And more generally, we have:
\begin{theorem}[Harper~\cite{harpertypicalchar}, 2023]\label{thmweightlow}
In the above setting with $r$ a large prime, uniformly for any $1 \leq x \leq r$, any $0 \leq q \leq 1$, and any multiplicative function $h(n)$ that has absolute value 1 on primes and absolute value at most 1 on prime powers, we have
$$ \frac{1}{r-1} \sum_{\chi \; \text{mod} \; r} |\sum_{n \leq x} h(n) \chi(n)|^{2q} \ll \Biggl(\frac{x}{1 + (1-q)\sqrt{\log\log(10L)}} \Biggr)^q , $$
where $L = L_r := \min\{x,r/x\}$.

Uniformly for any $1 \leq x \leq T$ (with $T \in \R$ large), any $0 \leq q \leq 1$, and any multiplicative function $h(n)$ that has absolute value 1 on primes and absolute value at most 1 on prime powers, we have
$$ \frac{1}{T} \int_{0}^{T} |\sum_{n \leq x} h(n) n^{it}|^{2q} dt \ll \Biggl(\frac{x}{1 + (1-q)\sqrt{\log\log(10L_T)}} \Biggr)^q  , $$
where $L_T := \min\{x,T/x\}$.
\end{theorem}
Recall here that a function $h$ is multiplicative if $h(mn) = h(m)h(n)$ for all coprime $m,n$ (total multiplicativity required this for {\em all} $m,n$). In particular, the M\"obius function $\mu(n)$ is multiplicative, and so a permissible choice of $h(n)$.

As a consequence of all this, provided that $L_r \rightarrow \infty$ as $r \rightarrow \infty$ (which means $x \rightarrow \infty$ but $x=o(r)$) we get $\frac{1}{r-1} \sum_{\chi \; \text{mod} \; r} |\sum_{n \leq x} \chi(n)| \ll \frac{\sqrt{x}}{(\log\log(10L_r))^{1/4}} = o(\sqrt{x})$, and so {\em we must have $|\sum_{n \leq x} \chi(n)| = o(\sqrt{x})$ for most $\chi$}. The same holds in the case of zeta sums $\sum_{n \leq x} n^{-it}$, provided that $L_T \rightarrow \infty$ as $T \rightarrow \infty$. And likewise for the sums with multiplicative weights.

\vspace{12pt}
In view of the ``Fourier flip'', the appearance of $L_r$ and $L_T$ in Theorem \ref{thmplainlow} is very natural, since up to scaling the sums up to $x$ there are essentially the same (in a loose sense) as the corresponding sums up to $r/x$ and $T/x$. When we discuss the proofs in section \ref{secmainntproofs} we will see how $L_r$ and $L_T$ arise, and it isn't too hard to show that if $0.01r \leq x \leq 0.99r$ (say) then indeed $\frac{1}{r-1} \sum_{\chi \; \text{mod} \; r} |\sum_{n \leq x} \chi(n)|^{2q} \asymp x^q$, with no extra cancellation. See the introduction of the original paper \cite{harpertypicalchar} for more discussion of this. In the more general setting of Theorem \ref{thmweightlow} there will not typically be any Fourier flip, and one might hope (for many $h$) that the bounds should be valid in a stronger form with $\log\log(10L)$ replaced by $\log\log(10x)$. We will revisit this point later.

\vspace{12pt}
The proofs of these theorems involve two types of arguments. The first, which is where the better than squareroot cancellation is really produced, shows that one has the same type of behaviour in an appropriate random model setting, of so-called {\em random multiplicative functions}. We will discuss this a little in section \ref{secrmf}, but it is not our focus here. See the author's paper~\cite{harperrmf3} for a recent survey of this random side of the theory. The other type of argument, which we will concentrate on, is of a more analytic and number theoretic flavour (although heavily motivated by the probabilistic side). This works to show, in a somewhat general way, that one can {\em transfer the results for averages of random multiplicative functions to deterministic families of multiplicative functions, like Dirichlet characters or $n^{-it}$}. We will discuss this in section \ref{secmainntproofs}.

\subsection{Potential application: the M\"obius function in short intervals and arithmetic progressions}
When Wintner~\cite{wintner} initiated his study of Rademacher random multiplicative functions (see section \ref{secrmf}), a key motivation was to model and gain insight into the behaviour of a specific multiplicative function, namely the M\"obius function $\mu(n)$. Nowadays it seems fairly safe to say that, except perhaps by ``coincidence'', the sequence of values of a Rademacher random multiplicative function do {\em not} generally provide a good model for $\mu(n)$. The M\"obius function is too special, with various properties flowing from its connection with Riemann zeta zeros that do not have analogues on the random multiplicative side. {\em But this doesn't mean that the random multiplicative model has nothing to say about $\mu(n)$.} As we will describe now, better than squareroot cancellation ideas originating in random multiplicative functions can not only predict, but (potentially) rigorously prove, results about $\mu(n)$ that are beyond the reach of other techniques. The key point is that these ideas are not applied directly to the single function $\mu(n)$, but to other {\em families of multiplicative functions} (like Dirichlet characters or $n^{-it}$ twisted by $\mu(n)$).

\vspace{12pt}
One of the most important questions in analytic number theory is to understand the distribution of primes in short intervals, or in arithmetic progressions. A close cousin of this question is to do the same thing for the values of the M\"{o}bius function. It seems very safe to conjecture that for any fixed $\epsilon > 0$, as $x \rightarrow \infty$ one should have
$$ \sum_{\substack{x < p \leq x+x^{\epsilon}, \\ p \; \text{prime}}} 1 = (1+o(1)) \int_{x}^{x+x^{\epsilon}} \frac{dt}{\log t} , $$
and
$$ \sum_{\substack{p \leq x, \\ p \equiv a \; \text{mod} \; r}} 1 = (1+o(1)) \frac{1}{\phi(r)} \int_{2}^{x} \frac{dt}{\log t}  , \; r \leq x^{1-\epsilon} , \; (a,r)=1 ; $$
and likewise
$$ \sum_{x < n \leq x+x^{\epsilon}} \mu(n) = o(x^{\epsilon}) , \;\;\; \text{and} \;\;\; \sum_{\substack{n \leq x, \\ n \equiv a \; \text{mod} \; r}} \mu(n) = o(x/r) , \; r \leq x^{1-\epsilon} . $$

Assuming the Riemann Hypothesis or the Generalised Riemann Hypothesis, these estimates are known to hold provided $\epsilon > 1/2$ (in fact one can count primes in intervals of length $\gg \sqrt{x} \log^{2}x$, and work of Maier and Montgomery~\cite{maiermont} obtains M\"obius cancellation in intervals of length $\sqrt{x} \log^{A}x$ for a certain constant $A$). For the M\"obius function, a consequence of the breakthrough work of Matom\"{a}ki and Radziwi\l\l~\cite{matradz} is the existence of positive and negative M\"{o}bius values (or, strictly speaking, values of the closely related Liouville function) in intervals of length $C\sqrt{x}$, for a large constant $C$. But we cannot handle intervals of length $\sqrt{x}$, or smaller. This is fundamentally due to squareroot cancellation, and not better than that, in the proofs: they discard the oscillations of certain terms and end up relying on only the horizontal, as opposed to vertical, distribution of zeros of $L$-functions; or (essentially equivalently) end up relying on estimates for absolute values of Dirichlet polynomials, rather than the polynomials themselves as complex numbers; and this seems to limit the maximum amount of cancellation available. In particular, Legendre's famous conjecture that there should be a prime between any two squares seems out of reach, since this would be equivalent (for $x$ large) to finding primes in intervals of length $(\sqrt{x}+1)^2 - (\sqrt{x})^2 \sim 2\sqrt{x}$.

It turns out that a certain conjectural extension of Theorem \ref{thmweightlow} would allow one to estimate sums of $\mu(n)$ in intervals of length (just) $o(\sqrt{x})$, or in arithmetic progressions mod $r$ with $r/\sqrt{x}$ (just) tending to infinity. In particular, {\em this would resolve the M\"obius function analogue of Legendre's conjecture}. The required statement is:
\begin{conjecture}[Harper~\cite{harpertypicalchar}, 2023]\label{conjlongchar}
For all $0 \leq q \leq 1$ and any fixed $A > 0$, we should have (for $r$ a large prime)
$$ \frac{1}{r-1} \sum_{\chi \; \text{mod} \; r} |\sum_{n \leq x} \mu(n) \chi(n)|^{2q} \ll_{A} (\frac{x}{1 + (1-q)\sqrt{\log\log x}} )^q \;\;\; \forall \; x \leq r^A , $$
and (for $T \in \R$ large)
$$ \frac{1}{2T} \int_{-T}^{T} |\sum_{n \leq x} \mu(n) n^{it}|^{2q} dt \ll_{A} (\frac{x}{1 + (1-q)\sqrt{\log\log x}} )^q \;\;\; \forall \; x \leq T^A . $$
\end{conjecture}
In addition to the extended ranges of $x$, we emphasise that in the denominators here we have $\log\log x$ rather than $\log\log L_r$ or $\log\log L_T$.

In fact, for the applications discussed we would ``only'' need Conjecture \ref{conjlongchar} when $A \approx 2$ (note that Theorem \ref{thmweightlow} already proves it when $A < 1$). Taking the short interval problem as an example, Perron's formula (multiplicative Fourier inversion) {\em on the $0$-line} implies that
$$ \left| \sum_{x < n \leq x+y} \mu(n) \right| \approx \left| \frac{1}{2\pi i} \int_{-i(x/y)}^{i(x/y)} (\sum_{n \leq 2x} \frac{\mu(n)}{n^s}) \frac{x^{s} ((1+y/x)^s - 1)}{s} ds \right| \lesssim \frac{y}{x} \int_{-x/y}^{x/y} |\sum_{n \leq 2x} \frac{\mu(n)}{n^{it}}| dt . $$
Provided that $y \leq x^{0.51}$, say, we can apply Conjecture \ref{conjlongchar} with $T = x/y \geq x^{0.49}$ and $A = 1/0.49$ to obtain (something like)
$$ \left| \sum_{x < n \leq x+y} \mu(n) \right| \ll \frac{\sqrt{x}}{(\log\log x)^{1/4}} . $$
This bound is $o(y)$, as desired, provided that $y$ grows faster than $\frac{\sqrt{x}}{(\log\log x)^{1/4}}$.

\vspace{12pt}
Conjecture \ref{conjlongchar} remains open, apart from the piece already proved in Theorem \ref{thmweightlow}, and is no doubt itself very challenging. Indeed, the full statement would imply the Riemann Hypothesis, since taking $y=x^{1-\epsilon}$ and $A = 1/\epsilon$ we could deduce roughly squareroot cancellation for partial sums of $\mu(n)$ (which is well known to be equivalent to the Riemann Hypothesis). However, this is not obviously the case when $A \approx 2$. And in any event it would be extremely interesting to prove the conjecture, with the consequences for $\mu(n)$ in short intervals and arithmetic progressions, even if one needed to assume the (Generalised) Riemann Hypothesis to do so. Note that the conjecture is still a statement about absolute values of Dirichlet polynomials and character sums, not about extra cancellation coming from their phases. In support of the conjecture, we have:
\begin{theorem}[Wang and Xu~\cite{wangxu}, 2025]\label{thmlongmob}
Let $r$ be a large prime. Assume the truth of the Generalised Riemann Hypothesis for all Dirichlet $L$-functions modulo $r$, and of a suitable form of the Ratios Conjecture for all primitive Dirichlet $L$-functions modulo $r$. Then for all $0 \leq q \leq 1$ and any fixed $A > 0$, we have
$$ \frac{1}{r-1} \sum_{\chi \; \text{mod} \; r} |\sum_{n \leq x} \lambda(n) \chi(n)|^{2q} \ll_{A} (\frac{x}{1 + (1-q)\sqrt{\log\log x}} )^q \;\;\; \forall \; x \leq r^A , $$
where $\lambda(n)$ is the Liouville function.
\end{theorem}
Note that the Liouville function is just the {\em totally} multiplicative variant of the M\"obius function, in other words $\lambda(n) := (-1)^{\Omega(n)}$ for {\em all} $n$ (without any condition of being squarefree), where $\Omega(n)$ denotes the total number of prime factors (with multiplicity) of $n$. Wang and Xu ultimately work with $\lambda(n)$ for technical convenience, but their arguments should certainly go through for $\mu(n)$ as well (and they already work out many steps of their proof explicitly in that case).

We will say a little about the proof of this theorem in section \ref{secwangxu}. We do not attempt to precisely state the Ratios Conjecture, which asserts that averages of ratios of (Dirichlet) $L$-functions should behave in the way suggested by assuming their zeros obey random matrix theory heuristics. This is unlikely to be proved any time soon, but it seems of some significance that Conjecture \ref{conjlongchar}, whose origins have no apparent connection with random matrix theory, is also consistent with those ideas. As Wang and Xu~\cite{wangxu} remark, ``Morally, [our] paper passes from random matrices to random multiplicative functions''.

\section{Aside: Other cases of better than squareroot cancellation}
Given the title of this survey, it seems remiss not to pause and briefly note some other number theoretic and analytic instances of better than squareroot cancellation. As the reader will see, although the phenomenon occurs in various places and may be exploited to great effect there, most of these are of a quite different nature to what happens in Theorems \ref{thmplainlow} and \ref{thmweightlow}. Firstly, in these examples one often sees not just better than squareroot cancellation, but {\em much} better than squareroot cancellation. More importantly, the cause of the behaviour is very different. The extra cancellation is not produced by multiplicativity or any subtle connections with randomness, but by the presence of quite rigid structure.

We can consider the additive character sum (or Dirichlet kernel) $\sum_{n \leq x} e(n\theta)$, where $e(\cdot) := e^{2\pi i \cdot}$ is the complex exponential. The trivial pointwise bound for this is $x$, and its second moment may easily be computed (using Parseval's identity/orthogonality of additive characters) to be $\int_{0}^{1} |\sum_{n \leq x} e(n\theta)|^2 d\theta = \lfloor x \rfloor$. So, by any reckoning, squareroot cancellation would suggest that for most $0 \leq \theta \leq 1$ we should have $|\sum_{n \leq x} e(n\theta)| \asymp \sqrt{x}$. But this sum is simply a geometric progression, and may be computed explicitly for any given $\theta$. We find that $|\sum_{n \leq x} e(n\theta)| = |\frac{e(\lfloor x \rfloor\theta) - 1}{e(\theta) - 1}| \ll \min\{x, \frac{1}{\Vert \theta \Vert}\}$, where $\Vert \theta \Vert$ is the distance from $\theta$ to the nearest integer. In particular, whenever $1/x^{1/10} \leq \theta \leq 1 - 1/x^{1/10}$, say (which is certainly ``most'' $\theta$), we get $|\sum_{n \leq x} e(n\theta)| \ll x^{1/10}$. We can also note that the first moment $\int_{0}^{1} |\sum_{n \leq x} e(n\theta)| d\theta \ll \int_{0}^{1} \min\{x, \frac{1}{\Vert \theta \Vert}\} d\theta \ll \log x$, again much smaller than squareroot size. The cause of this is the extreme additive structure of the set $\{1,2,...,\lfloor x \rfloor \}$ over which we sum, interacting with the additive nature of the $e(n\theta)$. The enormously better than squareroot cancellation that one gets here drives many common techniques in analytic number theory, such as the {\em completion of sums} method.

A related example arises from additive character sums over sets defined by digit constraints. For example, suppose that $b \in \N\backslash\{1\}$ is some base; that $(d_k)_{k \in \mathcal{K}}$ is some prescribed vector of base $b$ digits (so $\mathcal{K}$ is a small subset of $\{0,1,...,K-1\}$, and each of the $d_k$ lies between $0$ and $b-1$); and $\mathcal{D} := \{n = \sum_{k=0}^{K-1} \alpha_k b^k : 0 \leq \alpha_k \leq b-1, \; \alpha_k = d_k \; \forall \; k \in \mathcal{K}\}$. Then $\mathcal{D}$ is a subset of $\{0,1,...,b^K - 1\}$, with cardinality $\#\mathcal{D} = b^{K-\#\mathcal{K}}$, and squareroot cancellation would suggest that for most $\theta$ we should have $|\sum_{n \in \mathcal{D}} e(n\theta)| \asymp \sqrt{\#\mathcal{D}}$ (and then also $\int_{0}^{1} |\sum_{n \in \mathcal{D}} e(n\theta)| d\theta \asymp \sqrt{\#\mathcal{D}}$). But Bourgain~\cite{bourgainprimes} (for $b=2$, and Swaenepoel~\cite{swaenepoel} for general base $b$) has shown that the sums are typically much smaller, indeed $\int_{0}^{1} |\sum_{n \in \mathcal{D}} e(n\theta)| d\theta \leq \frac{(\#\mathcal{D})^{1+o(1)}}{b^K}$, where the $o(1)$ term tends to zero with the proportion of prescribed digits. Actually this follows, with some careful work, from the preceding example of the Dirichlet kernel, since $|\sum_{n \in \mathcal{D}} e(n\theta)| = \prod_{\substack{k=0, \\ k \notin \mathcal{K}}}^{K-1} |\sum_{n=0}^{b-1} e(nb^{k}\theta)|$. The very strong cancellation here makes it feasible to apply the Hardy--Littlewood circle method to attack two variable additive problems involving $\mathcal{D}$. For example, Bourgain~\cite{bourgainprimes} was able to count primes with a small positive proportion of binary digits fixed, by writing this count as $\int_{0}^{1} (\sum_{n \in \mathcal{D}} e(n\theta)) (\sum_{\substack{p \leq b^{K}, \\ p \; \text{prime}}} e(-p\theta)) d\theta$, and bounding the contribution from $\theta$ in the ``minor arcs'' $\mathfrak{m}$ by $\sup_{\theta \in \mathfrak{m}} |\sum_{\substack{p \leq b^{K}, \\ p \; \text{prime}}} e(-p\theta)| \cdot \int_{0}^{1} |\sum_{n \in \mathcal{D}} e(n\theta)| d\theta \leq \sup_{\theta \in \mathfrak{m}} |\sum_{\substack{p \leq b^{K}, \\ p \; \text{prime}}} e(-p\theta)| \cdot \frac{(\#\mathcal{D})^{1+o(1)}}{b^K}$. If one has a small power saving for $\sup_{\theta \in \mathfrak{m}} |\sum_{\substack{p \leq b^{K}, \\ p \; \text{prime}}} e(-p\theta)|$, which is possible with a good choice of $\mathfrak{m}$, then the minor arc contribution will be negligible compared with the anticipated size of the count (which is $\asymp \frac{\#\mathcal{D}}{\log(b^K)}$ under mild assumptions on the prescribed digits, the logarithm coming from the density of the primes). Notice that if $\int_{0}^{1} |\sum_{n \in \mathcal{D}} e(n\theta)| d\theta$ were anything close to $\sqrt{\#\mathcal{D}}$, as one would expect for a non-structured set, such an argument would completely fail. This is why one cannot apply this strategy to resolve e.g. the binary Goldbach problem.

Exponential sums over sets like the squarefree numbers or the $k$-free numbers typically enjoy much better than squareroot cancellation, because they can be obtained from the natural numbers by rather rigid sieve processes (with few sieving steps), and then inherit the better than squareroot cancellation of the Dirichlet kernel. More explicitly, we have $\sum_{\substack{n \leq x, \\ n \; \text{squarefree}}} e(n\theta) = \sum_{n \leq x} e(n\theta) \sum_{d: d^2 | n} \mu(d)$, with the M\"obius function playing its classical inclusion-exclusion role to detect the squarefree condition. Swapping the sums, we get $\sum_{d \leq \sqrt{x}} \mu(d) \sum_{n \leq x : d^2 | n} e(n\theta) = \sum_{d \leq \sqrt{x}} \mu(d) \sum_{m \leq x/d^2} e(md^{2}\theta)$, and using Parseval's identity we find for example that $\int_{0}^{1} |\sum_{D < d \leq \sqrt{x}} \mu(d) \sum_{m \leq x/d^2} e(md^{2}\theta)|^2 d\theta = \sum_{n \leq x} (\sum_{\substack{ d > D, \\ d^2 | n}} \mu(d))^2 \leq \sum_{n \leq x} d(n) \sum_{\substack{ d > D, \\ d^2 | n}} 1 \ll \frac{x \log^{O(1)}x}{D}$. Here $d(n) = \sum_{d|n} 1$ was the divisor function. So for most $\theta$, e.g. for a proportion 0.99 of $\theta$, the contribution from $d > D$ is $\ll \sqrt{\frac{x \log^{O(1)}x}{D}}$. In particular, this is significantly smaller than $\sqrt{x}$ if we take $D = x^{0.1}$, say, and we are left with only the short sum of Dirichlet kernels $\sum_{d \leq x^{0.1}} \mu(d) \sum_{m \leq x/d^2} e(md^{2}\theta)$. By elaboration of such ideas, Balog and Ruzsa~\cite{balogruzsa} showed that the first moment $\int_{0}^{1} |\sum_{\substack{n \leq x, \\ n \; \text{squarefree}}} e(n\theta)| d\theta \asymp x^{1/3}$ (and so the typical size of this exponential sum is $\ll x^{1/3}$). For $k$-free numbers with $k \geq 3$ the first moment is even smaller (of order $x^{1/(k+1)}$), reflecting the increasingly structured nature of these sets.

As a final and slightly different example, we mention the small size of the error terms in lattice point counting problems, such as the classical Dirichlet divisor problem (counting points under a hyperbola) and Gauss circle problem (counting points inside a circle). See e.g. the recent paper~\cite{lamzouridivisor} of Lamzouri for a discussion of such issues.

\section{Random multiplicative functions}\label{secrmf}
Let $(f(p))_{p \; \text{prime}}$ be independent Steinhaus random variables, i.e. independent random variables distributed uniformly on the complex unit circle $\{|z|=1\}$. We define a {\em Steinhaus random multiplicative function} $f : \N \rightarrow \C$, by setting $f(n) := \prod_{p^{a} || n} f(p)^{a}$ for all natural numbers $n$ (where $p^a || n$ means that $p^a$ is the highest power of the prime $p$ that divides $n$, so $n = \prod_{p^{a} || n} p^{a}$). Thus $f$ is a random function taking values in the complex unit circle, that is totally multiplicative, i.e. satisfies $f(nm) = f(n)f(m)$ for all $n,m$. This is the type of random multiplicative function most relevant to this survey.

An alternative model, and the one originally introduced by Wintner~\cite{wintner}, is to let $(f(p))_{p \; \text{prime}}$ be independent Rademacher random variables (i.e. taking values $\pm 1$ with probability $1/2$ each), and then set $f(n) := \prod_{p |n} f(p)$ for all squarefree $n$, and $f(n) = 0$ when $n$ is not squarefree. This gives a {\em Rademacher random multiplicative function}. As discussed previously, Wintner intended this as a direct model for (at least some features of) the M\"obius function $\mu(n)$, although nowadays one perhaps doesn't think of it quite like that. Note that both functions are supported on squarefree numbers, and both are multiplicative (but not totally multiplicative),  i.e. satisfy $f(nm) = f(n)f(m)$ and $\mu(nm) = \mu(n)\mu(m)$ for all {\em coprime} $n,m$.

\vspace{12pt}
By construction, a Steinhaus random multiplicative function shares the total multiplicativity property of Dirichlet characters $\chi(n)$ and continuous characters $n \mapsto n^{-it}$, and all of these take values in the complex unit circle (except when $n$ and $r$ are not coprime in the case of $\chi$). They also enjoy similar orthogonality properties. Now one takes an expectation $\E$ in place of averaging over $\chi$ or $t$, and we may calculate that
$$ \E f(n) \overline{f(m)} = \E \left(\prod_{p^{a} || n} f(p)^{a} \right) \left(\prod_{p^{a} || m} f(p)^{-a} \right) = \E \prod_{p| mn} f(p)^{a(n,p) - a(m,p)} = \prod_{p| mn} \E f(p)^{a(n,p) - a(m,p)} , $$
where $a(n,p), a(m,p)$ are the exponents of $p$ in the prime factorisations of $n,m$. Since $f(p)^{a(n,p) - a(m,p)}$ is uniformly distributed on the unit circle except when $a(n,p) = a(m,p)$, we get $\E f(n) \overline{f(m)} = \textbf{1}_{a(n,p)=a(m,p) \; \forall \, p} = \textbf{1}_{n=m}$. Notice this holds for {\em all} $n,m$, without any of the size restrictions we had earlier for $\chi$ and $n^{-it}$. This is a crucial point that we shall return to later. A related observation is that a Steinhaus random multiplicative function will (almost surely) {\em not} exhibit any of the periodicity or symmetry properties of $\chi(n)$ or $n^{-it}$. Indeed this would not make sense, since for a random multiplicative function there is no attached ``conductor'' $r$ or $T$.

\vspace{12pt}
Since Steinhaus random multiplicative functions share two of the three key properties of Dirichlet characters and continuous characters (i.e. multiplicativity and some form of orthogonality, but not periodicity), one may hope that results about Steinhaus random multiplicative functions could be a useful guide to what happens for $\chi(n)$ and $n^{-it}$. {\em It turns out this is often the case, provided the question one studies only involves values of $n$ (on the deterministic side) not ``too close'' to the conductor $r$ or $T$, so one doesn't ``see'' the period too much\footnote{For a couple of examples, see the papers of Granville and Soundararajan~\cite{gransoundlcs} and Lamzouri~\cite{lamzouri2dzeta}. The same point is briefly discussed in Sarnak's letter~\cite{sarnakpd} on positive definite $L$-functions.}.}

Although random multiplicative functions share several properties of number theoretic functions, making them rather interesting, they also enjoy a crucial property that the number theoretic functions lack--- the exact independence of the underlying $f(p)$. This makes them much more susceptible to analysis with probabilistic tools and ideas. In the author's opinion, the general scheme of work in this area is to first prove results for $f(n)$ using any available tools, in particular leveraging the independence of the $(f(p))_{p \; \text{prime}}$, and then see how far these can be transported to deterministic results for $\chi(n)$ or $n^{-it}$ (or potentially for other families of multiplicative functions as well).

\vspace{12pt}
Thanks to orthogonality, it is easy to see that
$$ \E|\sum_{n \leq x} f(n)|^{2} = \sum_{n,m \leq x} \textbf{1}_{n=m} = \lfloor x \rfloor $$
for a Steinhaus random multiplicative function $f(n)$. Squareroot cancellation would then suggest that we should typically (e.g. with probability 0.99) have $|\sum_{n \leq x} f(n)| \asymp \sqrt{x}$. But a result of the author~\cite{harperrmflowmoments} shows that we get more cancellation than this. Indeed, uniformly for all large $x$ and all real $0 \leq q \leq 1$ (possibly depending on $x$) we have
\begin{equation}\label{rmforder}
\E|\sum_{n \leq x} f(n)|^{2q} \asymp \left( \frac{x}{1 + (1-q)\sqrt{\log\log x}} \right)^q ,
\end{equation}
with an extra saving factor of $\log\log x$ in the denominator. Notice the strong (and of course not coincidental!) similarity of this statement with Theorems \ref{thmplainlow} and \ref{thmweightlow}. {\em On the random side one has an order of magnitude result rather than (as yet) just an upper bound on the number theoretic side, suggesting that the apparently strange shape of the bounds in Theorems \ref{thmplainlow} and \ref{thmweightlow} may be the correct answer.} As one should expect since there is no notion of conductor $r$ or $T$ here, the denominator in \eqref{rmforder} is $\log\log x$ rather than $\log\log L$.

We will not say too much about the proof of \eqref{rmforder}, referring the reader to e.g. Harper~\cite{harperrmf3} for a survey of this and related results. The introduction of the original paper~\cite{harperrmflowmoments} also contains a (hopefully) fairly accessible discussion, although lacking some more recent technical simplifications. But one key point is that the proof splits quite neatly into two, somewhat uneven, parts. In the first part, one shows that for a suitable parameter $P \leq x$ (which can be chosen with a lot of flexibility, at least for the upper bound argument) we have
\begin{equation}\label{eqnboundbyprod}
\E |\sum_{n \leq x} f(n)|^{2q} \ll x^{q} \E\Biggl(\frac{1}{\log P} \int_{-1/2}^{1/2} |F_{P}^{\text{rand}}(1/2+it)|^2 dt \Biggr)^{q} + \text{small error} ,
\end{equation}
where $F_{P}^{\text{rand}}(s) := \prod_{p \leq P} (1 - \frac{f(p)}{p^s})^{-1}$ is the random Euler product corresponding to $f(n)$. Then one shows that
\begin{equation}\label{eqnboundforprod}
\E (\frac{1}{\log P} \int_{-1/2}^{1/2} |F_{P}^{\text{rand}}(1/2+it)|^2 dt )^{q} \approx (\frac{1}{1 + (1-q)\sqrt{\log\log P}} )^q .
\end{equation}
Provided $P$ may be chosen so that $\log\log P \gg \log\log x$, the combination of these estimates gives the desired upper bound for $\E|\sum_{n \leq x} f(n)|^{2q}$.

The second step of the argument, estimating $\E (\frac{1}{\log P} \int_{-1/2}^{1/2} |F_{P}^{\text{rand}}(1/2+it)|^2 dt )^{q}$, is certainly the more challenging part, and involves probabilistic ideas related to random walk and critical multiplicative chaos. Very loosely, the additional saving factor $\log\log P$ reflects non-trivial correlations between the Euler products $F_{P}^{\text{rand}}(1/2+it), F_{P}^{\text{rand}}(1/2+it')$ for nearby $t, t'$. But when it comes to transferring the result to the number theoretic side, it turns out that we can insert this estimate for $\E (\frac{1}{\log P} \int_{-1/2}^{1/2} |F_{P}^{\text{rand}}(1/2+it)|^2 dt )^{q}$ entirely as a black box. In particular, we emphasise that the argument to bound $\frac{1}{r-1} \sum_{\chi \; \text{mod} \; r} |\sum_{n \leq x} \chi(n)|^{2q}$ or $\frac{1}{r-1} \sum_{\chi \; \text{mod} \; r} |\sum_{n \leq x} \mu(n) \chi(n)|^{2q}$, or the analogous $t$ averages, {\em directly exploits an estimate for the random quantity $\E (\frac{1}{\log P} \int_{-1/2}^{1/2} |F_{P}^{\text{rand}}(1/2+it)|^2 dt )^{q}$}.

To make a transfer to the deterministic side, the new work is to establish a deterministic analogue of \eqref{eqnboundbyprod} with e.g. $\frac{1}{r-1} \sum_{\chi \; \text{mod} \; r} |\sum_{n \leq x} \chi(n)|^{2q}$ on the left hand side, {\em but still with a random quantity like $\E (\frac{1}{\log P} \int_{-1/2}^{1/2} |F_{P}^{\text{rand}}(1/2+it)|^2 dt )^{q}$ on the right}. In fact, we end up producing a slightly weaker and more complicated variant of \eqref{eqnboundbyprod}, where the right hand side can still be bounded by inserting our estimate for $\E (\frac{1}{\log P} \int_{-1/2}^{1/2} |F_{P}^{\text{rand}}(1/2+it)|^2 dt )^{q}$. This is the argument that we shall try to explain in the next section.

\section{Proof ideas for Theorems \ref{thmplainlow} and \ref{thmweightlow}}\label{secmainntproofs}
We try to indicate how one can obtain a deterministic version of \eqref{eqnboundbyprod}. For concreteness, we shall mostly write things in the particular case of $\frac{1}{r-1} \sum_{\chi \; \text{mod} \; r} |\sum_{n \leq x} \chi(n)|^{2q}$, and to simplify the writing we shall hereafter write $\Echar$ to denote averaging over all Dirichlet characters mod $r$. (Thus if $W(\chi)$ is a function, then $\Echar W := \frac{1}{r-1} \sum_{\chi \; \text{mod} \; r} W(\chi)$.) See also the introduction of the original paper of Harper~\cite{harpertypicalchar}, for a similar but more compressed and technical discussion.

\vspace{12pt}
To set the scene, we first observe that if $p_1, ..., p_k, p_{k+1}, ..., p_l$ are any (not necessarily distinct) primes such that $\prod_{j=1}^{k} p_j, \prod_{j=k+1}^{l} p_j < r$, we have an equality
\begin{eqnarray}\label{charrmfexp}
\Echar \prod_{j=1}^{k} \chi(p_j) \prod_{j=k+1}^{l} \overline{\chi(p_j)} & = & \Echar \chi(\prod_{j=1}^{k} p_j) \overline{\chi(\prod_{j=k+1}^{l} p_j)} = \textbf{1}_{\prod_{j=1}^{k} p_j \equiv \prod_{j=k+1}^{l} p_j \; \text{mod} \; r} \nonumber \\
& = & \textbf{1}_{\prod_{j=1}^{k} p_j = \prod_{j=k+1}^{l} p_j} = \E \prod_{j=1}^{k} f(p_j) \prod_{j=k+1}^{l} \overline{f(p_j)} .
\end{eqnarray}
Here we used multiplicativity and orthogonality \eqref{eqncharorth} of Dirichlet characters; the size restrictions on $\prod_{j=1}^{k} p_j$ and $\prod_{j=k+1}^{l} p_j$ (to crucially replace a congruence mod $r$ by an equality); and finally the orthogonality property $\E f(n) \overline{f(m)} = \textbf{1}_{n=m}$ of Steinhaus random multiplicative functions. Notice we are exploiting precisely the properties\footnote{We can also emphasise here the reason for assuming the conductor $r$ to be prime: if it were not, then \eqref{charrmfexp} can fail for any primes $p_j$ dividing $r$, where the left hand side would vanish but the right hand side need not. However, it should be straightforward to extend everything to the case of more general $r$ (at least those not having ``too many'' small prime divisors, e.g. assuming that the Euler totient value $\phi(r) \gg r$), by modifying the definition of $f(n)$ on primes dividing $r$ to preserve \eqref{charrmfexp}. If one looks at $t$ averages rather than character averages, this sort of issue does not arise.} shared by $\chi$ and the random multiplicative model $f$.

More generally, \eqref{charrmfexp} implies that the $\Echar$ average of polynomials in the $\chi(p_j)$ will exactly match the $\E$ average of the same polynomials in $f(p_j)$, provided the ``degree'' (or length, or complexity) of all the terms remains under control, i.e. provided we have $\prod_{j=1}^{k} p_j, \prod_{j=k+1}^{l} p_j < r$ in all the arising terms. Since one can approximate quite general functions by polynomials, this gives us some hope of being able to match up the distribution of terms involving $\chi$ with the corresponding terms involving $f$. To understand the behaviour of $\sum_{n \leq x} \chi(n)$, it is clearly the terms $(\chi(p))_{p \leq x}$ that need to be analysed. Unfortunately, unless $x$ is very small compared with $r$ (e.g. a small power of $\log r$) any direct argument along these lines, trying to match the (joint) distribution of the $\chi(p)$ with that of the $f(p)$, will not succeed. This is because the degree of the polynomials needed to produce good approximations grows quickly, thus the condition that $\prod_{j=1}^{k} p_j, \prod_{j=k+1}^{l} p_j < r$ will quickly be violated as $x$ (and so the number and size of the primes $p_j$) grows. As we will explain shortly, there is also a clear conceptual reason why a full comparison of joint distributions must fail unless $x$ is small.

But the equality \eqref{charrmfexp} will be fundamental to our more subtle argument, that can be made to work and deliver strong bounds even for larger $x$ (as it is fundamental, whether used implicitly or explicitly, in essentially all arguments involving averages of character sums). At a high level, the key observation is that {\em one does not need a good understanding of the joint distribution of all the $(\chi(p))_{p \leq x}$ in order to understand $\Echar |\sum_{n \leq x} \chi(n)|^{2q}$}. The different $\chi(p)$ are ``blended together'' inside the sum, moreover we are helped because we are only seeking understanding at the level of moments (albeit fractional moments), and only seeking upper bounds.

\vspace{12pt}
We now change direction for a moment, and discuss how \eqref{eqnboundbyprod} can be obtained for random multiplicative functions. As we shall see, this again exploits multiplicativity and orthogonality (which would be available for $\chi$ as well, at least to some extent), but also crucially uses a {\em conditioning} procedure that does not extend immediately (unless $x$, or really the parameter $P$, is very small). This part of our discussion will have a bit more probabilistic flavour than the rest of this survey.

Let $P(n)$ denote the largest prime factor of $n$. Provided our parameter $P$ isn't too close to $x$, the contribution to $\sum_{n \leq x} f(n)$ from $P$-smooth numbers $n$ (i.e. those $n$ with $P(n) \leq P$) may be discarded into the error term in \eqref{eqnboundbyprod}, because the second moment of this contribution is fairly small (there aren't many $P$-smooth numbers). This leaves us to deal with $\E|\sum_{n \leq x, P(n) > P} f(n)|^{2q}$. Each number in this sum factors uniquely in the form $nm$, where $n > 1$ has all of its prime factors $> P$ (and therefore actually $n > P$), and $m$ is $P$-smooth. Then using multiplicativity of $f$, we can write
$$ \sum_{n \leq x, P(n) > P} f(n) = \sum_{\substack{P < n \leq x, \\ p|n \Rightarrow p > P}} f(n) \sum_{\substack{m \leq x/n, \\ P(m) \leq P}} f(m) . $$
{\em Note that the inner sums are entirely determined by the $(f(p))_{p \leq P}$}. Thus if we let $\tilde{\E}$ denote expectation conditional on the values $(f(p))_{p \leq P}$ (i.e. expectation with those values treated as fixed and the $(f(p))_{p > P}$ remaining random, so the conditional expectation of any quantity is a function of the values $(f(p))_{p \leq P}$), we observe
\begin{eqnarray}
\tilde{\E} \Biggl| \sum_{n \leq x, P(n) > P} f(n) \Biggr|^{2} & = & \sum_{\substack{P < n_1, n_2 \leq x, \\ p|n_{1}n_{2} \Rightarrow p > P}} \tilde{\E}\Biggl( f(n_1) \Bigl(\sum_{\substack{m_1 \leq x/n_1, \\ P(m_1) \leq P}} f(m_1) \Bigr) \overline{f(n_2)} \Bigl(\sum_{\substack{m_2 \leq x/n_2, \\ P(m_2) \leq P}} \overline{f(m_2)} \Bigr) \Biggr) \nonumber \\
& = & \sum_{\substack{P < n \leq x, \\ p|n \Rightarrow p > P}} \Biggl| \sum_{\substack{m \leq x/n, \\ P(m) \leq P}} f(m) \Biggr|^2 , \nonumber
\end{eqnarray}
because the $(f(p))_{p > P}$ {\em remain independent with mean zero}, and so the terms $f(n_1), f(n_2)$ {\em remain orthogonal under the conditional expectation $\tilde{\E}$}. This step crucially uses the independence of the $f(p)$, more specifically the fact that fixing all of the $(f(p))_{p \leq P}$ changes nothing about the distribution of the $(f(p))_{p > P}$. As we shall emphasise shortly, such full independence cannot possibly extend to the $\chi(p)$ unless $P$ is very small.

To finish the derivation of \eqref{eqnboundbyprod}, by the Tower Property of conditional expectation (which here is simply Fubini's theorem, breaking up the multiple ``integration'' $\E$ into separate integrations corresponding to the $(f(p))_{p \leq P}$ and the $(f(p))_{p > P}$) we can write $\E|\sum_{n \leq x, P(n) > P} f(n) |^{2q} = \E \tilde{\E} | \sum_{n \leq x, P(n) > P} f(n) |^{2q}$. So applying H\"{o}lder's inequality to the {\em conditional expectation} $\tilde{\E}$ only, we get
$$ \E\Biggl|\sum_{n \leq x, P(n) > P} f(n) \Biggr|^{2q} \leq \E \Biggl( \tilde{\E} \Biggl| \sum_{n \leq x, P(n) > P} f(n) \Biggr|^{2} \Biggr)^{q} = \E \Biggl( \sum_{\substack{P < n \leq x, \\ p|n \Rightarrow p > P}} \Biggl| \sum_{\substack{m \leq x/n, \\ P(m) \leq P}} f(m) \Biggr|^2 \Biggr)^{q} . $$
Having reached this point, one can perform some smoothing of the $n$ sum and use a suitable form of Parseval's identity to relate the right hand side to an Euler product average like $\E (\frac{1}{\log P} \int_{-1/2}^{1/2} |F_{P}^{\text{rand}}(1/2+it)|^2 dt )^{q}$. The factor $\frac{1}{\log P}$ reflects the density of $P$-rough numbers $n$ (i.e. numbers $n$ with all their prime factors $> P$) on the interval $(P,x]$.

We mentioned a couple of times that we cannot hope to match up the joint distribution of the $(\chi(p))_{p \leq P}$ or $(\chi(p))_{p \leq x}$ with the corresponding $f(p)$, or perform the sort of conditioning that would follow from this, unless $P$ or $x$ is very small. This is simply because {\em there are not enough characters $\chi$ mod $r$} for the distributions to possibly remain close. Indeed, the smallest possible non-zero ``probability'' of any ``event'' involving the $\chi$ will be $1/(r-1)$, if precisely one character $\chi$ satisfies the given conditions. In contrast, on the random side one has a much richer product probability space. As a simple but very relevant example, on the random side we have $\p(\Re f(p) \geq 0 \; \forall \; p \leq P) = (1/2)^{\#\{p \leq P\}} \approx (1/2)^{P/\log P}$, which will be much smaller than $1/(r-1)$ (yet still non-zero) as soon as $P \gg (\log r) \log\log r$. So if we try to condition (freeze) the values $\chi(p)$ for too many primes $p$, even at the crude level of the signs of $\Re \chi(p)$, we cannot hope for conditioning on such configurations to match the calculations in the random multiplicative case. The failure of \eqref{charrmfexp} when the condition $\prod_{j=1}^{k} p_j, \prod_{j=k+1}^{l} p_j < r$ is violated is a manifestation of the failure of such matching, and of the lack of genuine independence of the values $\chi(p)$.

\vspace{12pt}
To rescue the situation, we observe that we did not really need to condition on the exact values of all $(f(p))_{p \leq P}$ on the random multiplicative side. Again, these values are ``blended together'' inside $\sum_{n \leq x, P(n) > P} f(n)$, which in any case we are only trying to understand at the level of upper bounds for fractional moments. Looking at the final statement of \eqref{eqnboundbyprod} gives us a good idea of the amount of information we really needed to condition on (freeze) in $\tilde{\E}$ to carry through the proof. Since we ultimately bound $\sum_{\substack{P < n \leq x, \\ p|n \Rightarrow p > P}} | \sum_{\substack{m \leq x/n, \\ P(m) \leq P}} f(m) |^2$ by something like $\frac{1}{\log P} \int_{-1/2}^{1/2} |F_{P}^{\text{rand}}(1/2+it)|^2 dt$, the proof could proceed similarly if we conditioned on any (potentially much lesser) information about the $(f(p))_{p \leq P}$ that roughly fixes the value of $\int_{-1/2}^{1/2} |F_{P}^{\text{rand}}(1/2+it)|^2 dt$. Putting this another way, we want to ultimately have a bound roughly like $\E (\frac{1}{\log P} \int_{-1/2}^{1/2} |F_{P}^{\text{rand}}(1/2+it)|^2 dt )^{q}$, and this will be unchanged provided enough randomness remains in the outer averaging $\E$ to roughly preserve the distribution of $\frac{1}{\log P} \int_{-1/2}^{1/2} |F_{P}^{\text{rand}}(1/2+it)|^2 dt$. This should certainly be possible with much less kept on the outside than a full averaging over each individual $(f(p))_{p \leq P}$.

It turns out that $|F_{P}^{\text{rand}}(1/2+it)| \approx \exp\{\Re \sum_{p \leq P} \frac{f(p)}{p^{1/2+it}}\}$ does not usually change much when $t$ varies by less than about $1/\log P$. So if we replace $\tilde{\E}$ in the above description by a {\em coarser conditioning}, only on the approximate values of $\Re \sum_{p \leq P} \frac{f(p)}{p^{1/2+it}}$ at a net of $t$ values with spacing roughly $1/\log P$, this should fix enough information on the inside of the $q$-th power (and therefore leave enough averaging on the outside of the $q$-th power) that we still end up with an overall bound of roughly the desired shape $\E (\frac{1}{\log P} \int_{-1/2}^{1/2} |F_{P}^{\text{rand}}(1/2+it)|^2 dt )^{q}$.

Let us try to make this even more explicit on the probabilistic side, to make the passage back to the character sum problem more obvious. If $\tilde{\E}$ is now {\em any} conditioning that still {\em only involves the values $(f(p))_{p \leq P}$}, (so the values $f(n)$ on $P$-rough numbers $n$ remain orthogonal under $\tilde{\E}$), then exactly the same Tower Property argument as above would yield a bound $\E|\sum_{n \leq x, P(n) > P} f(n) |^{2q} \leq \E ( \tilde{\E} \sum_{\substack{P < n \leq x, \\ p|n \Rightarrow p > P}} | \sum_{\substack{m \leq x/n, \\ P(m) \leq P}} f(m) |^2 )^{q}$. Notice that $\tilde{\E}$ remains on the inside at present, since if we do not condition on all the exact values of $(f(p))_{p \leq P}$ then $\tilde{\E} \sum_{\substack{P < n \leq x, \\ p|n \Rightarrow p > P}} | \sum_{\substack{m \leq x/n, \\ P(m) \leq P}} f(m) |^2$ may no longer be the same as $\sum_{\substack{P < n \leq x, \\ p|n \Rightarrow p > P}} | \sum_{\substack{m \leq x/n, \\ P(m) \leq P}} f(m) |^2$. But by the same smoothing and Parseval argument as before, which is rather general and makes little use of the nature of $\tilde{\E}$, we can upper bound $\sum_{\substack{P < n \leq x, \\ p|n \Rightarrow p > P}} | \sum_{\substack{m \leq x/n, \\ P(m) \leq P}} f(m) |^2$ here (up to acceptable average errors) by something like $\frac{1}{\log P} \int_{-1/2}^{1/2} |F_{P}^{\text{rand}}(1/2+it)|^2 dt$. Thus we end up with $\E ( \tilde{\E} \frac{1}{\log P} \int_{-1/2}^{1/2} |F_{P}^{\text{rand}}(1/2+it)|^2 dt )^{q}$. Provided $\tilde{\E}$ conditions on enough information to roughly fix the value of $\frac{1}{\log P} \int_{-1/2}^{1/2} |F_{P}^{\text{rand}}(1/2+it)|^2 dt$, for example it conditions on the approximate values of $\Re \sum_{p \leq P} \frac{f(p)}{p^{1/2+it}}$ at a net of $t$ values with spacing $1/\log P$, then $\tilde{\E}$ now disappears and we again arrive at an upper bound like the desired $\E (\frac{1}{\log P} \int_{-1/2}^{1/2} |F_{P}^{\text{rand}}(1/2+it)|^2 dt )^{q}$.

\vspace{12pt}
Now let us return to the deterministic setting, and more explicitly number theoretic considerations. For character sums, our version of ``conditioning'' will simply be breaking up $\sum_{\chi \; \text{mod} \; r}$ according to the behaviour of certain quantities, in this case something like the values of $\Re \sum_{p \leq P} \frac{\chi(p)}{p^{1/2+it}}$ at all points $t$ in a suitable net. The earlier discussion about the coarseness of this conditioning roughly translates into the fact that, at least quite often, we want the subsums into which $\sum_{\chi \; \text{mod} \; r}$ is broken to still contain a large number of characters. At first pass one could think of this breaking up as simply splitting up the characters into classes using indicator functions, but in fact (precisely so that we can exploit the match between averages of polynomials in $\chi(p)$ and in $f(p)$, as in \eqref{charrmfexp}) we will do the splitting using smooth functions.

Let us give a little more detail, referring to the original paper of Harper~\cite{harpertypicalchar} for much more precise discussion. We write
$$ \Echar |\sum_{n \leq x} \chi(n)|^{2q} = \sum_{\textbf{j}} \Echar G_{\textbf{j}}(\chi) |\sum_{n \leq x} \chi(n) |^{2q} , $$
where $G_{\textbf{j}}(\chi)$ are smooth, non-negative functions satisfying $\sum_{\textbf{j}} G_{\textbf{j}}(\chi) \equiv 1$ that approximately pick out all characters $\chi$ for which sums like $\Re \sum_{p \leq P} \frac{\chi(p)}{p^{1/2+it}}$, at a suitable net of $t$, have a given collection of values. (In fact, for the later steps one should first remove those $n$ from the sum that don't have a large prime factor, which can be done with exactly the same second moment argument as in the random case.) Here the sum $\sum_{\textbf{j}}$ plays roughly the same role as the outer averaging $\E$ in the equality $\E|\sum_{n \leq x, P(n) > P} f(n) |^{2q} = \E \tilde{\E} | \sum_{n \leq x, P(n) > P} f(n) |^{2q}$, and the $\Echar G_{\textbf{j}}(\chi)$ play the same role as the conditional expectations $\tilde{\E}$. To make the analogy even closer, we can rewrite our equality slightly as
$$ \Echar |\sum_{n \leq x} \chi(n)|^{2q} = \sum_{\textbf{j}} \sigma(\textbf{j}) \E^{\textbf{j}} \Biggl|\sum_{n \leq x} \chi(n) \Biggr|^{2q} , $$
where for each $\textbf{j}$ in the outer sum we set $\sigma(\textbf{j}) := \Echar G_{\textbf{j}}(\chi)$, and (for all functions $W(\chi)$) set $\E^{\textbf{j}} W := \frac{1}{\sigma(\textbf{j})} \Echar G_{\textbf{j}}(\chi) W(\chi)$. These manipulations ensure that the weighted outer summation $\sum_{\textbf{j}} \sigma(\textbf{j})$ is actually an averaging (i.e. $\sum_{\textbf{j}} \sigma(\textbf{j}) = 1$), likewise each of the inner ``conditional'' $\E^{\textbf{j}}$ is actually an averaging (i.e. $\E^{\textbf{j}} 1 = 1$).

We can apply H\"{o}lder's inequality to each inner average $\E^{\textbf{j}} |\sum_{n \leq x} \chi(n) |^{2q}$, exactly as we did with $\tilde{\E} | \sum_{n \leq x, P(n) > P} f(n) |^{2q}$, and find
$$ \Echar |\sum_{n \leq x} \chi(n)|^{2q} \leq \sum_{\textbf{j}} \sigma(\textbf{j}) \Biggl( \E^{\textbf{j}} \Biggl|\sum_{n \leq x} \chi(n) \Biggr|^{2} \Biggr)^q . $$
This leaves us to understand quantities like $\E^{\textbf{j}} |\sum_{n \leq x} \chi(n) |^{2} = \frac{1}{\sigma(\textbf{j})} \Echar G_{\textbf{j}}(\chi) |\sum_{n \leq x} \chi(n) |^{2}$. Since the $G_{\textbf{j}}$ are nice smooth functions, we can then approximate each $G_{\textbf{j}}(\chi)$ by a polynomial in the $\chi(p)$ (up to acceptable errors), and also expand out the square $|\sum_{n \leq x} \chi(n) |^{2}$. Thus we can (potentially) invoke \eqref{charrmfexp} to replace $\sigma(\textbf{j})$ by $\sigma^{\text{rand}}(\textbf{j}) := \E G_{\textbf{j}}(f)$, and replace $\Echar G_{\textbf{j}}(\chi) |\sum_{n \leq x} \chi(n) |^{2}$ by $\E G_{\textbf{j}}(f) |\sum_{n \leq x} f(n) |^{2}$. This leads to a bound of the shape
\begin{equation}\label{firstderand}
\Echar |\sum_{n \leq x} \chi(n)|^{2q} \lesssim \sum_{\textbf{j}} \sigma^{\text{rand}}(\textbf{j}) \Biggl( \E^{\textbf{j}, \text{rand}} \left|\sum_{n \leq x} f(n) \right|^{2} \Biggr)^{q} + \text{small error} ,
\end{equation}
where we wrote $\E^{\textbf{j}, \text{rand}} W(f)$ to denote $\frac{1}{\sigma^{\text{rand}}(\textbf{j})} \E G_{\textbf{j}}(f) W(f)$. {\em Notice this is precisely the type of bound that we wanted, with the deterministic quantity $\Echar |\sum_{n \leq x} \chi(n)|^{2q}$ on the left, and everything on the right hand side genuinely random.} If $G_{\textbf{j}}(\chi)$ (and thus $G_{\textbf{j}}(f)$) is constructed to only involve the values $(\chi(p))_{p \leq P}$, then inside $\E^{\textbf{j}, \text{rand}} \left|\sum_{n \leq x} f(n) \right|^{2}$ the values of $f$ on $P$-rough numbers will remain orthogonal, and so $\left|\sum_{n \leq x} f(n) \right|^{2}$ can be replaced by something like $\sum_{\substack{P < n \leq x, \\ p|n \Rightarrow p > P}} | \sum_{\substack{m \leq x/n, \\ P(m) \leq P}} f(m) |^2$. Then we can invoke our usual smoothing and Parseval type argument, and end up with a bound like
\begin{equation}\label{secondderand}
\Echar |\sum_{n \leq x} \chi(n)|^{2q} \lesssim x^q \sum_{\textbf{j}} \sigma^{\text{rand}}(\textbf{j}) \Biggl( \E^{\textbf{j}, \text{rand}} \frac{1}{\log P} \int_{-1/2}^{1/2} |F_{P}^{\text{rand}}(1/2+it)|^2 dt \Biggr)^{q} .
\end{equation}

\vspace{12pt}
We end with a few more remarks about the derivation and exploitation of \eqref{secondderand}.

Most importantly, we remarked above that we could {\em potentially} invoke \eqref{charrmfexp} to replace $\sigma(\textbf{j})$ by $\sigma^{\text{rand}}(\textbf{j})$, and replace $\Echar G_{\textbf{j}}(\chi) |\sum_{n \leq x} \chi(n) |^{2}$ by $\E G_{\textbf{j}}(f) |\sum_{n \leq x} f(n) |^{2}$. In fact, we can do this precisely when {\em the total lengths of the involved character sums are $< r$}. Since $G_{\textbf{j}}$ is a smooth function that is supposed to (approximately) detect the joint values of $\Re \sum_{p \leq P} \frac{\chi(p)}{p^{1/2+it}}$ at $\approx \log P$ points $t$, it turns out that $G_{\textbf{j}}(\chi)$ can be acceptably approximated by a character sum of length $e^{\log^{O(1)}P}$. Meanwhile, $|\sum_{n \leq x} \chi(n) |^{2}$ contributes a character sum (and its conjugate) of length $x$. So we can make our argument work provided $x e^{\log^{O(1)}P} < r$, allowing us to choose $P \approx \exp\{\log^{c}L_r\}$ (recall that $L_r = \min\{x,r/x\}$). We then see\footnote{If we had not switched to ``conditioning'' on the sums $\Re \sum_{p \leq P} \frac{\chi(p)}{p^{1/2+it}}$, and persevered with a direct approach of using $G_{\textbf{j}}(\chi)$ to roughly detect the individual values of all $(\chi(p))_{p \leq P}$, then the approximating character sum would need to have length more like $e^P$. Thus we would need a condition like $x e^{P} < r$, and could not generally take $P$ large enough to conclude a sharp bound in Theorem \ref{thmplainlow}.} that we will ultimately obtain a saving $\log\log P \asymp \log\log L_r$, as claimed in Theorem \ref{thmplainlow}. As already discussed in the Introduction, and see also section \ref{secwangxu} below, in the case of $\sum_{n \leq x} \chi(n)$ (but perhaps not for some weighted sums) this appearance of $L_r$ in our bounds reflects a real feature of the problem, namely the ``Fourier flip'' symmetry of the sums. We see again that the length restriction in \eqref{charrmfexp} reflects genuine aspects of the behaviour of these characters, whereby they do not entirely mimic the pure random multiplicative model as the period $r$ comes into play.

Note that in the derivation of \eqref{secondderand}, we pass to the random side very early. Indeed this is already done in \eqref{firstderand}, after some basic manipulations to set up our ``conditioning'' and apply H\"older's inequality, and then a matching of character sum averages with random multiplicative function expectations. In particular, all of the smoothing and Parseval steps to reach an Euler product, and then all work with Euler products, is done on the random side. Organising the argument in this way is not essential, one could work directly with Dirichlet characters for much longer and only match up with the random side at the end (to allow a final invocation of \eqref{eqnboundforprod}). However, although it may make things a little less familiar for a number theory audience, there seem to be advantages to proceeding as we do. For example, working with characters one could never successfully analyse Euler product averages (unless the Euler products were extremely short). Thus it would be necessary to work with shorter character sums that approximate the Euler product, in a suitable sense, instead. This style of argument is now well established in analytic number theory, especially the study of moments of $L$-functions, and could be implemented--- but at the price of additional technicality. By passing quickly to random multiplicative functions, one quickly gets access to perfect orthogonality (with no attached size restrictions on the variables) and independence, and can handle objects like full Euler products directly.

As written, the right hand side of \eqref{secondderand} does not look exactly like the right hand side $x^{q} \E(\frac{1}{\log P} \int_{-1/2}^{1/2} |F_{P}^{\text{rand}}(1/2+it)|^2 dt )^{q}$ of \eqref{eqnboundbyprod}. The point is that, provided the functions $G_{\textbf{j}}$ and the parameter $P$ are set properly (to localise all the sums $\Re \sum_{p \leq P} \frac{f(p)}{p^{1/2+it}}$ precisely enough, to some values depending only on $\textbf{j}$), then the values of $\frac{1}{\log P} \int_{-1/2}^{1/2} |F_{P}^{\text{rand}}(1/2+it)|^2 dt$ that contribute to $\E^{\textbf{j}, \text{rand}} \frac{1}{\log P} \int_{-1/2}^{1/2} |F_{P}^{\text{rand}}(1/2+it)|^2 dt = \frac{1}{\sigma^{\text{rand}}(\textbf{j})} \E G_{\textbf{j}}(f) \frac{1}{\log P} \int_{-1/2}^{1/2} |F_{P}^{\text{rand}}(1/2+it)|^2 dt$ will all be more or less the same (depending only on $\textbf{j}$). Then in \eqref{secondderand} we get
\begin{eqnarray}
&& \sum_{\textbf{j}} \sigma^{\text{rand}}(\textbf{j}) \Biggl( \E^{\textbf{j}, \text{rand}} \frac{1}{\log P} \int_{-\frac{1}{2}}^{\frac{1}{2}} |F_{P}^{\text{rand}}(\frac{1}{2}+it)|^2 \Biggr)^{q} \nonumber \\
& = & \sum_{\textbf{j}} \E G_{\textbf{j}}(f) \Biggl( \E^{\textbf{j}, \text{rand}} \frac{1}{\log P} \int_{-\frac{1}{2}}^{\frac{1}{2}} |F_{P}^{\text{rand}}(\frac{1}{2}+it)|^2 \Biggr)^{q} \nonumber \\
& \approx & \sum_{\textbf{j}} \E G_{\textbf{j}}(f) \Biggl( \frac{1}{\log P} \int_{-\frac{1}{2}}^{\frac{1}{2}} |F_{P}^{\text{rand}}(\frac{1}{2}+it)|^2 \Biggr)^{q} , \nonumber
\end{eqnarray}
where the first equality comes from unpacking the definition of $\sigma^{\text{rand}}(\textbf{j})$, and the second (approximate) equality comes from the fact that $\frac{1}{\log P} \int_{-1/2}^{1/2} |F_{P}^{\text{rand}}(1/2+it)|^2 dt$ is more or less the same wherever $G_{\textbf{j}}(f)$ is supported\footnote{In the rigorous proofs of Theorems \ref{thmplainlow} and \ref{thmweightlow}, this is the step with the most demanding conditions on the functions $G_{\textbf{j}}$.} (or, in reality, wherever it is not very small). Then recalling that $\sum_{\textbf{j}} G_{\textbf{j}}(f) \equiv 1$, we recover a bound like $\E(\frac{1}{\log P} \int_{-1/2}^{1/2} |F_{P}^{\text{rand}}(1/2+it)|^2 dt )^{q}$.

\vspace{12pt}
Finally, we take this opportunity to flag up a couple of aspects of the proof that a reader seeking to generalise to other situations should especially keep in mind. Most importantly, when constructing the functions $G_{\textbf{j}}(\chi)$ for conditioning, one must ensure that they fix enough information (i.e. approximately localise enough sums like $\Re \sum_{p \leq P} \frac{\chi(p)}{p^{1/2+it}}$) so that at the end of the proof (when one has passed to the random side), whatever quantity appears in place of $\frac{1}{\log P} \int_{-1/2}^{1/2} |F_{P}^{\text{rand}}(1/2+it)|^2 dt$ is approximately fixed. For example, in the proof of Theorem \ref{thmplainlow} it is not actually $\frac{1}{\log P} \int_{-1/2}^{1/2} |F_{P}^{\text{rand}}(1/2+it)|^2 dt$ that appears (we wrote this as a slight expository simplification), one rather gets something like $\frac{1}{\log P} \int_{-\infty}^{\infty} \frac{|F_{P}^{\text{rand}}(1/2+it)|^2}{|1/2 + it|^2} dt$. The part of this integral where $|t| > \log^{0.01}P$, say, can fairly easily be discarded thanks to the decay from $\frac{1}{|1/2 + it|^2}$, but we cannot simply discard all $|t| > 1/2$. Thus in $G_{\textbf{j}}(\chi)$ we actually ``condition'' on $\Re \sum_{p \leq P} \frac{\chi(p)}{p^{1/2+it}}$ (or later $\Re \sum_{p \leq P} \frac{f(p)}{p^{1/2+it}}$) at a net of $t$ values covering all $|t| \leq \log^{0.01}P$, not just $|t| \leq 1/2$. In other situations, different types of integral will arise, and the important sums and range of $t$ may be different again. We emphasise in particular that, although the {\em distribution} of $\frac{1}{\log P} \int_{n-1/2}^{n+1/2} |F_{P}^{\text{rand}}(1/2+it)|^2 dt$ is exactly the same for any shift $n$, the actual values of these integrals (for any given realisation of the random multiplicative function $f(n)$) will change as $n$ varies. Thus it is necessary to include all $|t| \leq \log^{0.01}P$ (or whatever the relevant range of $t$ proves to be) in the conditioning construction.

We also note that in the setting of Theorem \ref{thmweightlow}, if the twist function $h(n)$ is not totally multiplicative (only multiplicative) then some small adjustments are required when taking account of the prime square contribution in the conditioning process. This is worked out explicitly by Harper~\cite{harpertypicalchar}. Similar considerations may be relevant if one tried to adapt the entire set-up to situations where total multiplicativity or perfect orthogonality breaks down, for example replacing the average $\frac{1}{r-1} \sum_{\chi \; \text{mod} \; r}$ by an average over quadratic characters, or families arising from higher degree $L$-functions.

\section{The work of Wang and Xu}\label{secwangxu}
Recall that the key source of the restrictions $x \leq r$ and $x \leq T$ in our unconditional main theorems was the need, in the proofs, to always be working with (squares of) character sums or Dirichlet polynomials of length $\leq r$ and $\leq T$, so that their averages over $\chi$ or $t$ would match the corresponding averages involving random multiplicative functions $f$.

We should also note again that, in general, this is not simply an artefact of the proof, but may reflect real changes in behaviour. For example, if one studies the moments $\frac{1}{r-1} \sum_{\chi \; \text{mod} \; r} |\sum_{n \leq x} \chi(n)|^{2q}$ of unweighted character sums with $x$ large enough compared with $r$, then the principal character $\chi_0 = \textbf{1}_{(n,r)=1}$ alone gives a contribution $\gg x^{2q}/r$, which may be much larger than $x^q$. This is a somewhat artificial objection, which one might look to resolve simply by excluding the principal character from the sum (see e.g. the work of Szab\'o~\cite {szaboupper,szabolower} on high moments of character sums, assuming the Generalised Riemann Hypothesis for the upper bounds). However, if one does this it still makes little sense to consider these unweighted sums with $x > r$, since the sum of $\chi(n)$ over a complete period $r$ is exactly zero for any $\chi \neq \chi_0$, so all such sums anyway reduce to being of length $< r$. Recall also that for e.g. $\frac{1}{r-1} \sum_{\chi \; \text{mod} \; r} |\sum_{n \leq x} \mu(n) \chi(n)|^{2q}$, one hopes to obtain a saving $\log\log x$ in the denominator in Conjecture \ref{conjlongchar}, as opposed to the $\log\log(10L) = \log\log(10\min\{x,r/x\})$ that must arise for $\sum_{n \leq x} \chi(n)$. Thus, to obtain strong bounds for $\frac{1}{r-1} \sum_{\chi \; \text{mod} \; r} |\sum_{n \leq x} \mu(n) \chi(n)|^{2q}$ where $x$ is permitted to approach or exceed $r$, one should expect to input some extra information about $\mu(n)$ (which in particular distinguishes it from the constant weight 1).

\vspace{12pt}
In proving Theorem \ref{thmlongmob}, Wang and Xu~\cite{wangxu} achieve exactly this on assuming the Generalised Riemann Hypothesis and the Ratios Conjecture. They work to show that whenever $x \leq r^A$ (with $A$ fixed but potentially large), and $G(\chi)$ is a short character sum (arising as a polynomial approximation to a function $G_{\textbf{j}}(\chi)$, as in section \ref{secmainntproofs}), one has
\begin{equation}\label{eqorthsub}
\frac{1}{r-1} \sum_{\chi \; \text{mod} \; r} G(\chi) |\sum_{n \leq x} \mu(n) \chi(n)|^2 = \E G(f) |\sum_{n \leq x} \mu(n) f(n)|^2 + \text{small error} .
\end{equation}
Given such a statement, (which would hold with no error term at all if $x < r$ and $G(\chi)$ has length $< r/x$, by orthogonality as in \eqref{charrmfexp}), one can otherwise follow the proof of Theorem \ref{thmweightlow} more or less exactly to obtain Theorem \ref{thmlongmob} (although with some non-trivial calculation to check that the combination of all the small, but non-zero, error terms remains under control). Note that the reason one can have a denominator $\log\log x$ for $\frac{1}{r-1} \sum_{\chi \; \text{mod} \; r} |\sum_{n \leq x} \mu(n) \chi(n)|^{2q}$ is because $G(\chi)$ will always be allowed to have length up to a small power of $x$ in \eqref{eqorthsub}, rather than $\min\{x,r/x\}$ when invoking orthogonality.

\vspace{12pt}
To obtain \eqref{eqorthsub}, most of Wang and Xu's\footnote{As discussed in the Introduction, Wang and Xu's~\cite{wangxu} final result (our Theorem \ref{thmlongmob}) handles $\sum_{n \leq x} \lambda(n) \chi(n)$ with $\lambda(n)$ the Liouville function, but the distinction between this and the M\"obius function is purely technical and is unimportant in almost all steps of their argument. In particular, they prove \eqref{eqorthsub} for both $\mu(n)$ and $\lambda(n)$, assuming the Generalised Riemann Hypothesis and the Ratios Conjecture.} work is with $\sum_{n \leq x} \mu(n) \chi(n)$. The short character sum $G(\chi)$ is simply broken up and treated termwise (so they actually work to show that $\frac{1}{r-1} \sum_{\chi \; \text{mod} \; r} \chi(n) \overline{\chi(m)} |\sum_{n \leq x} \mu(n) \chi(n)|^2 = \E f(n) \overline{f(m)} |\sum_{n \leq x} \mu(n) f(n)|^2 + \text{small error}$, uniformly for all $n,m$ up to a small power of $x$). Using Perron's formula, and assuming the Generalised Riemann Hypothesis to shift the line of integration close to 1/2, we have $\sum_{n \leq x} \mu(n) \chi(n) \approx \frac{1}{2\pi i} \int_{1/2 + \epsilon - iT}^{1/2 + \epsilon + iT} \frac{1}{L(s,\chi)} \frac{x^s}{s} ds$, where $T$ is suitably large (in terms of $r,x$) and $\epsilon > 0$ is a small constant. Expanding out, this implies that
\begin{eqnarray}
&& \frac{1}{r-1} \sum_{\chi \; \text{mod} \; r} \chi(n) \overline{\chi(m)} |\sum_{n \leq x} \mu(n) \chi(n)|^2 \nonumber \\
& \approx & \frac{1}{r-1} \sum_{\chi \; \text{mod} \; r} \frac{1}{(2\pi)^2} \int_{1/2 + \epsilon - iT}^{1/2 + \epsilon + iT} \int_{1/2 + \epsilon - iT}^{1/2 + \epsilon + iT} \frac{\chi(n) \overline{\chi(m)}}{L(s_1,\chi) \overline{L(s_2,\chi)}} \frac{x^{s_1}}{s_1} \frac{x^{\overline{s_2}}}{\overline{s_2}} ds_1 \overline{ds_2} . \nonumber
\end{eqnarray}
To understand this, it would clearly suffice to understand $\frac{1}{r-1} \sum_{\chi \; \text{mod} \; r} \frac{\chi(n) \overline{\chi(m)}}{L(s_1,\chi) \overline{L(s_2,\chi)}} = \frac{1}{r-1} \sum_{\chi \; \text{mod} \; r} \frac{\chi(n) \overline{\chi(m)}}{L(s_1,\chi) L(\overline{s_2},\overline{\chi})}$, and this is what can be done with some work assuming the Ratios Conjecture. (Actually one should exclude the principal character from the sum at this point, and reinsert it later.) Indeed, the terms $\chi(n), \overline{\chi(m)}$ are simply Dirichlet series coefficients of $L(z_1,\chi), L(z_2, \overline{\chi})$, so another application of Perron's formula or Mellin inversion implies that it would suffice to understand the behaviour of $\frac{1}{r-1} \sum_{\chi \; \text{mod} \; r} \frac{L(z_1,\chi) L(z_2, \overline{\chi})}{L(s_1,\chi) L(\overline{s_2},\overline{\chi})}$, with $z_1, z_2$ further complex variables. Wang and Xu~\cite{wangxu} assume the Ratios Conjecture in the form of a main term for this average (predicted by random matrix theory), and an error term giving a small fixed power saving in $r$, and which is reasonably uniform in $|z_1|, |z_2|, |s_1|, |s_2|$.

Note that the very rough size of the integrand above is $x^{1+2\epsilon}$. Here $x^{\epsilon} \leq r^{A\epsilon}$, and since $\epsilon$ may be taken arbitrarily small (thanks to the Generalised Riemann Hypothesis), the overall contribution from the Ratios Conjecture error term becomes $\ll_{A} \frac{x}{r^{\text{small power}}}$. This is acceptable for \eqref{eqorthsub}. The main term supplied by the Ratios Conjecture splits into two pieces. The first of these exactly corresponds to the random multiplicative function average (or diagonal contribution), as desired for \eqref{eqorthsub}. The second piece turns out to contribute negligibly because of the way it appears inside a contour integral, specifically because we remain strictly to the right of the 1/2 line in the variables $s_1, s_2$.

We emphasise some very particular properties of $\mu(n)$ that drove this argument, which would be unavailable for a general multiplicative weight $h(n)$. The most crucial feature is the expression for $\sum_{n \leq x} \mu(n) \chi(n)$ in terms of $\frac{1}{L(s,\chi)}$, where the $L$-function has nice properties (at least assuming strong conjectures) and can be understood. However, one could not do the same for all sums expressible in some way in terms of $L(s,\chi)$, as we noted earlier the analogue of Theorem \ref{thmlongmob} for the unweighted sum $\sum_{n \leq x} \chi(n)$ would be false in various ways. Thus the specific main term supplied by the Ratios Conjecture, and the nature of the contribution from the principal character (negligible when twisted by $\mu(n)$, but not when untwisted), are also both important.

\section{Final remarks}
We close this survey with some very brief further remarks.

The arguments leading to \eqref{rmforder} in the random case are quite robust, and have already been applied to various weighted sums $\sum_{n \leq x} f(n) a_n$ as well (and to other, less obviously related, situations). See e.g. the final section of the author's survey \cite{harperrmf3} for some examples. Given a result on the random side, the procedure from section \ref{secmainntproofs} for deducing corresponding deterministic results is also robust (although less so than the purely random steps), so one should be able to obtain better than squareroot cancellation bounds for character and zeta sums in various settings. One example that has already been worked out to some extent is that of short interval sums $\sum_{x < n \leq x+y} f(n)$ and $\sum_{x < n \leq x+y} \chi(n)$. Although the indicator function of a short interval is certainly not a multiplicative weight, one may think of it as somewhat close to that if $y$ isn't too far from $x$, and Caich~\cite{caichshort} has shown that better than squareroot cancellation persists for $\sum_{x < n \leq x+y} f(n)$ (i.e. the sum is typically $o(\sqrt{y})$) provided $\log(x/y)$ grows slower than $\sqrt{\log\log x}$. He obtains roughly analogous results for $\sum_{x < n \leq x+y} \chi(n)$, although now with somewhat complicated interaction between the sizes of $x,y$ and $r$. It would be interesting to explore this further, as well as other sums.

As already mentioned, at present Theorems \ref{thmplainlow} and \ref{thmweightlow} provide upper bounds only, whereas in the probabilistic situation \eqref{rmforder} we have matching upper and lower bounds. Motivated by this, the author~\cite{harpertypicalchar} conjectured that the bounds in Theorem \ref{thmplainlow} should be sharp, provided that $x \leq 0.99r$ in the character sum case (when $x$ is extremely close to $r$, the behaviour changes because of perfect periodicity). The best existing lower bounds for these moments differ from our upper bounds by powers of $\log L$, see e.g. section 1.5 of La Bret\`eche, Munsch and Tenenbaum~\cite{bretechemunschten}. (Note that we do {\em not} expect the bounds in Theorem \ref{thmweightlow} to be sharp in general, as for a typical multiplicative twist $h(n)$ one expects a saving $\log\log x$ rather than the weaker $\log\log L$.)

Matching lower bounds, if true, would have applications to non-vanishing of character sums, and to the closely related Dirichlet theta functions $\theta(s,\chi)$ at $s=1$, say. For example, if $\chi$ is an even Dirichlet character mod $r$ then $\theta(1, \chi) = \sum_{n=1}^{\infty} \chi(n) e^{-\pi n^{2}/r} \approx \sum_{n \leq \sqrt{r}} \chi(n)$, and it follows easily from Theorem \ref{thmplainlow} and partial summation that for any large prime $r$, and uniformly for $0 \leq q \leq 1$, we have
\begin{equation}\label{momenttheta}
\frac{1}{r-1} \sum_{\substack{\chi \; \text{mod} \; r, \\ \chi \; \text{even}}} |\theta(1, \chi)|^{2q} \ll \Biggl(\frac{\sqrt{r}}{1+(1-q)\sqrt{\log\log r}} \Biggr)^q .
\end{equation}
See Corollary 2 of Harper~\cite{harpertypicalchar} for a full proof of this. If we had sharp lower bounds in Theorem \ref{thmplainlow} when $x \approx \sqrt{r}$ (and when averaging only over even characters), this would essentially imply a matching lower bound in \eqref{momenttheta}. Comparing these estimates when e.g. $q=1/2$ and $q=2/3$ using H\"older's inequality, because the moments depend on $q$ only linearly in the exponent (provided $q$ is strictly away from 1) we could then deduce that $\theta(1,\chi) \neq 0$ for a positive proportion of (even) $\chi$ mod $r$. No such positive proportion non-vanishing result is currently known by any method. This lower bound question will be addressed in forthcoming work of the author.

The most tantalising prospect is perhaps to establish a result like Theorem \ref{thmlongmob} under less formidable hypotheses. This would be a considerable challenge, but it is worth noting that even if one follows Wang and Xu's~\cite{wangxu} strategy exactly there is no real need for a precise asymptotic for $\frac{1}{r-1} \sum_{\chi \; \text{mod} \; r} \frac{L(z_1,\chi) L(z_2, \overline{\chi})}{L(s_1,\chi) L(\overline{s_2},\overline{\chi})}$ for all $z_1,z_2,s_1,s_2$, since these ratios all appear with (some) averaging over those variables. One could also hope to access further averaging at the point where $L$-function computations are required, by suitably adjusting some earlier steps of the argument.

\vspace{12pt}
\noindent {\em Acknowledgements.} The author would like to thank Ben Green for his detailed comments and encouragement.

\end{document}